\documentstyle{amsppt}
\magnification1200 
\pagewidth{6.5 true in} 
\pageheight{9.25 true in}
\NoBlackBoxes

\topmatter
\title Large character sums: 
Pretentious characters and the Polya-Vinogradov theorem
\endtitle
\author Andrew Granville and K. Soundararajan
\endauthor
\rightheadtext 
{Pretentious characters and the Polya-Vinogradov theorem}
\address{D{\'e}partment  de Math{\'e}matiques et Statistique,
Universit{\'e} de Montr{\'e}al, CP 6128 succ Centre-Ville,
Montr{\'e}al, QC  H3C 3J7, Canada}\endaddress
\email{andrew{\@}dms.umontreal.ca}
\endemail
\address{Department of Mathematics, University of Michigan, Ann Arbor,
Michigan 48109, USA} \endaddress \email{ksound{\@}umich.edu}
\endemail
\thanks{Le premier auteur est partiellement soutenu par une bourse
de la Conseil  de recherches en sciences naturelles et eng\' enie
du Canada. The second  author is partially supported by the
National Science Foundation.}
\endthanks
\endtopmatter

\def\cbar{\overline{\chi}}

\document
\head 1. Introduction \endhead

\noindent 
The best bound known for character sums was given independently by
G. P{\' o}lya and I.M. Vinogradov in 1918 (see [4], p.135-137).  
For any non-principal Dirichlet character $\chi\pmod{q}$ we 
let $$M(\chi) := \max_{x} \Big|\sum_{n\le x} \chi(n)\Big|,$$ and 
then the P{\' o}lya-Vinogradov inequality reads 
$$
M(\chi) \ll \sqrt{q} \log q. \tag{1.1}
$$
There has been no subsequent improvement in this 
inequality other than in the implicit constant.  
Moreover it is believed that (1.1) will be 
difficult to improve since it is possible 
(though highly unlikely) that there is an infinite sequence of primes 
$q\equiv 1 \pmod 4$ for which $(\frac pq)=1$ for all $p<q^\epsilon$,
in which case $M((\frac{\cdot}{q})) \gg_{\epsilon} \sqrt{q} \log q$.

The unlikely possibility described above involves a quadratic character, 
and one might imagine that there are similar possibilities 
preventing one from improving (1.1) for higher order characters. 
Surprisingly, one of our main results shows that we can improve (1.1) 
for characters of odd, bounded order.

\proclaim{Theorem 1} If $\chi \pmod q$ is a 
primitive character of odd order $g$ then  
$$
M(\chi) \ll_g \sqrt{q} (\log q)^{1-\frac{\delta_g}{2} +o(1)},
$$
where $\delta_g= (1- \frac{g}{\pi} \sin \frac{\pi}{g})$.
\endproclaim  

Our proof of Theorem 1 is based on some technical 
results (described in the next section) which allow 
us to characterize characters $\chi$ for which $M(\chi)$ 
is large.  Our characterization reveals that 
there is a hidden structure among the characters 
having large $M(\chi)$.  One example of this 
structure is the following: 

\proclaim{Theorem 2} For $1\le j\le g$ let $\chi_j \pmod {q_j}$ be 
primitive characters (not necessarily distinct) with $q_j \le q$ for all $j$.  
We suppose that the product $\chi_1 \cdots \chi_g$ gives 
the principal character.  If $g$ is odd then we have that 
$$
\prod_{j=1}^{g} \frac{M(\chi_j)}{\sqrt{q_j}} \ll_g (\log q)^{g-\frac{1}{2g}}
+(\log q)^{g-\frac 17}.
$$
If $g$ is even then 
$$
\frac{M(\chi_g)}{\sqrt{q_g} \log q_g} +(\log q)^{-\frac {2(g-1)}7} 
\gg_g \prod_{j=1}^{g-1} 
\Big(\frac{M(\chi_j)}{\sqrt{q_j}\log q}\Big)^{2(g-1)}.
$$ 
\endproclaim 

Roughly speaking, the first part of Theorem 2 tells us that if 
$g$ is odd and $\chi_1 \cdots \chi_g =1$ then at least 
one of the $M(\chi_j)$ is small. In particular, taking 
$\chi_1=\dots=\chi_{k}=\chi$ with $k=g-1$, if 
$M(\chi)$ is large then $M(\chi^k)$ is small for 
small even integers $k$. The second part 
of Theorem 2 tells us that if $g$ is even 
and $M(\chi_1)$, $\ldots$, $M(\chi_{g-1})$ are all 
large, then so is $M(\chi_1 \cdots \chi_{g-1})$.  In particular,
taking $\chi_1=\dots=\chi_{k}=\chi$ with $k=g-1$, if 
$M(\chi)$ is large then $M(\chi^k)$ is also large for 
small odd numbers $k$.   

Another consequence of Theorems 1 and 2 
is that if $q$ is prime $\equiv 3 \pmod 
4$, and $M(\chi)\gg \sqrt{q}\log q$
for some character $\chi \pmod q$ of bounded order then
$M((\frac{\cdot}{q}))\gg \sqrt{q}\log q$
for the quadratic character $(\frac{\cdot}{q})$. 
One can deduce further results like this
from Theorem 2.

We give yet a third consequence.  
Suppose that $q_1$, $q_2$, $q_3$ are pairwise coprime, odd, 
squarefree integers in the interval $[Q,2Q]$, such that each 
$M((\frac{\cdot}{q_i}))\gg \sqrt{q_i}\log q_i$. Then 
$M((\frac{\cdot}{q_1q_2})) \ll \sqrt{q_1q_2}(\log (q_1q_2))^{6/7}$,
whereas $M((\frac{\cdot}{q_1q_2q_3})) \gg \sqrt{q_1q_2q_3}\log (q_1q_2q_3)$.
Similar results can be proved for products of four or more characters.

These bounds are larger than the expected maximal order of character 
sums.  In 1977 H.L. Montgomery and R.C. Vaughan [12]  
showed if the Generalized Riemann Hypothesis{\footnote {In their results 
and in ours (when indicated), the Riemann 
hypothesis for all Dirichlet $L$-functions is needed, not merely 
the Riemann hypothesis for the particular $L(s,\chi)$. } }
 is true 
then 
$$
M(\chi) \ll \sqrt{q} \log\log q. \tag{1.2}
$$
This bound is best possible, up to the evaluation of the constant, 
in view of R.E.A.C. Paley's 1932 result [13] that 
there are infinitely many positive integers $q$ such that  
$$
M((\tfrac {\cdot}{q})) \geq
\Big( \frac{e^\gamma}{\pi}+o(1)\Big) \sqrt{q} \log\log q, \tag{1.3}
$$
where $\gamma = 0.5772\ldots$ is the Euler-Mascheroni
constant.\footnote{Actually Paley's method gives the constant ``$1/2$'' not
``$e^\gamma$'', but such an improvement appeared subsequently in
several places, for example [1].}
Paley's result gives large character sums for a thin class of carefully 
constructed quadratic characters, and one may ask if for each large prime 
$q$ there are characters $\chi\pmod q$ with similarly large $M(\chi)$.  Our 
next result shows that there are indeed many such characters $\chi$, and 
moreover we can point these character sums in any given direction.

\proclaim{Theorem 3}  Let $q$ be a large prime and let 
$\theta \in (-\pi,\pi]$ be given.  There is an absolute 
constant $C_0$ such that for at least 
$q^{1-C_0/(\log \log q)^2}$ characters $\chi \pmod q$ with $\chi(-1)=-1$ 
we have 
$$
\sum_{n\le x}\chi(n) = e^{i\theta} \frac{e^{\gamma}}{\pi} \sqrt{q} \Big( 
\log \log q
+O((\log \log q)^{1/2})\Big)
$$
for all but $o(q)$ natural numbers $x\leq q$.
\endproclaim

In view of Theorem 3 it may be surprising that there are 
analogues of Theorems 1 and 2 which give a sharper upper bound 
than (1.2) for characters of small odd order. 

\proclaim{Theorem 4} Assume GRH.  If $\chi \pmod q$ is a 
primitive character of odd order $g$ then  
$$
M(\chi) \ll_g \sqrt{q} (\log\log q)^{1-\frac{\delta_g}{2} +o(1)},
$$
where $\delta_g= (1- \frac{g}{\pi} \sin \frac{\pi}{g})$.
\endproclaim 

On GRH, we can show that there exist arbitrarily 
large $q$ and primitive characters 
$\chi \pmod q$ of odd order $g$ such that
$$
M(\chi) \gg_g
\sqrt{q} (\log \log q)^{1-\delta_g -o(1)}. \tag{1.4} 
$$
We believe that the exponent $1-\delta_g =
\frac g{\pi} \sin \frac{\pi}{g}$ in (1.4) 
is best possible, and that Theorem 4 can be improved to 
attain this bound.  Perhaps this can be achieved by improving 
the bound given in Lemma 4.3.  It would also be interesting to obtain 
the lower bound (1.4) unconditionally.

\proclaim{Theorem 5} Assume GRH. For $1\le j\le g$ let $\chi_j 
\pmod {q_j}$ be 
primitive characters with $q_j\le q$ for all $j$. 
We suppose that the product $\chi_1 \cdots \chi_g$ gives 
the principal character.  If $g$ is odd then we have that 
$$
\prod_{j=1}^{g} \frac{M(\chi_j)}{\sqrt{q_j}} \ll_g 
(\log \log q)^{g-\frac{1}{2g}} + (\log \log q)^{g-\frac 17}.
$$
If $g$ is even then 
$$
\frac{M(\chi_g)}{\sqrt{q_g}\log \log q_g} + (\log \log q)^{-\frac {2(g-1)}7}
\gg_g \prod_{j=1}^{g-1}\Big(\frac{M(\chi_j)}{\sqrt{q_j} \log\log q} 
\Big)^{2(g-1)}.
$$
\endproclaim  

One can make deductions from Theorems 4 and 5 analogous to those consequences we 
gave after Theorems 1 and 2.

The known estimates  for character sums strongly resemble 
bounds for $L(1,\chi)$.  Unconditionally it is easy to show that 
$|L(1,\chi)|\ll \log q$.  Assuming GRH, J.E. Littlewood [11] 
proved that 
$$
L(1,\chi) \sim \prod_{p\le \log^2 q} \Big(1-\frac{\chi(p)}{p}\Big)^{-1}, 
\tag{1.5} 
$$
from which it follows easily that
$$
|L(1,\chi)\le (1+o(1)) 2 e^{\gamma}\log \log q. \tag{1.6}
$$
Apart from a factor of $2$, the bound (1.6) is best 
possible, since S.D. Chowla [3] showed that 
there exist arbitrarily large $q$ and characters $\chi \pmod q$ 
such that 
$$
|L(1,\chi)| \ge (1+o(1)) e^{\gamma}\log \log q.
$$
 
\proclaim{Theorem 6} Assume GRH.  If $\chi$ is 
a primitive character $\pmod q$ then 
$$ 
\Big|\sum_{n\le x}  \chi(n) \Big| 
\le \Big(\frac{2e^{\gamma}}{\pi} +o(1)\Big) \sqrt{q}\log \log q.
$$
Further 
$$
\Big| \sum_{x\le n\le x+y} \chi(n) \Big| \le 
\Big(\frac{4e^{\gamma}}{\pi \sqrt{3}} +o(1)\Big) \sqrt{q} \log \log q. 
$$
\endproclaim 

Regarding the second part of Theorem 6 we record that 
with minor modifications to the proof of Theorem 4 we 
may prove that for any angle $\theta \in (-\pi, \pi]$ and 
any large prime $q$ there are at least $q^{1-2/(\log \log q)^2}$ 
characters $\chi\pmod q$ with $\chi(-1)=1$ such that 
$$
\sum_{q/3 \le n \leq 2q/3} \chi(n) \sim e^{i\theta} 
\frac{2e^{\gamma}}{\pi \sqrt{3}} \sqrt{q}  
\log \log q. \tag{1.7}
$$
Theorem 6 places the situation for large character sums 
on the same footing as bounds for $L(1,\chi)$: the conditional 
$O$-results for character sums differ from Paley's $\Omega$ result 
by only a factor of $2$.  Moreover the maximal size of character 
sums in an interval $[x,x+y]$ is also determined up to a factor of $2$. 
It is believed that the $\Omega$ result represents the true 
extreme values of $L(1,\chi)$ (see [8] for arguments in the 
case when $\chi$ is quadratic).  Similarly we believe 
that (1.3) and (1.7) give the largest possible character sums.  

\proclaim{Conjecture 1} If $\chi$ is a primitive character 
$\pmod q$ then
$$
\Big| \sum_{n\leq x} \chi(n) \Big| \leq
\Big( \frac{e^\gamma}\pi +o(1)\Big)   \sqrt{q} \ \log\log q,  
$$
and 
$$
\Big|\sum_{x\le n\le x+y} \chi(n) \Big| 
\le \Big( \frac{2e^{\gamma}}{\pi \sqrt{3}} +o(1)\Big) 
\sqrt{q} \log \log q.
$$
\endproclaim

\head 2.  Detailed statement of results \endhead

\noindent If $\chi$ is a primitive character $\pmod q$ then the sum 
$\sum_{n\le x} \chi(n)$ has a Fourier expansion 
which is given quantitatively as (see P{\' o}lya [14]) 
$$
\sum_{n\leq x} \chi(n) = \frac{\tau (\chi)}{2i\pi} \sum\Sb n\in
{\Bbb Z} \\ 1\leq |n|\leq N \endSb \frac{\overline{\chi} (n)}{n}
\Big( 1- e(-\tfrac{nx}{q}) \Big) + O\Big(1+\frac{q\log q}N \Big).
$$
Here $\tau(\chi)$ is the usual Gauss sum (see section 4).  
Choosing $N=q$ above and noting that $L(1,\cbar)=\sum_{n\le q} 
\cbar(n)/n + O(1)$ we obtain that 
$$
\sum_{n\le x} \chi(n) = 
\frac{\tau(\chi)}{2\pi i} (1-\cbar(-1)) L(1,\cbar) 
- \frac{\tau(\chi)}{2\pi i} \sum_{n\le q} 
\frac{\cbar(n)}{n} \Big( e(-\tfrac {nx}{q})-\chi(-1)e(\tfrac{nx}{q})
\Big) + O(\sqrt{q}). \tag{2.1} 
$$   

All of our work here proceeds from the Fourier expansion (2.1).  
We wish to understand when the terms appearing in (2.1) can 
be large.  Littlewood's result (1.5) indicates that $L(1,\cbar)$ 
is large only when $\cbar(p) \approx 1$ for many small 
primes $p$.  We will find that the other terms appearing 
in (2.1) can be large only when $\chi(p)\approx \xi(p)$ for 
many small primes $p$, where $\xi$ is a character of small conductor.   
A. Hildebrand [10] first realized the possibility of such a 
result.

To formulate our results precisely we define for two 
characters $\chi$ and $\psi$  
$$
{\Bbb D}(\chi,\psi;y) := \Big( \sum_{p\le y} \frac{1-\text{Re }\chi(p) 
\overline{\psi}(p)}{p} \Big)^{\frac 12}. \tag{2.2} 
$$
We think of ${\Bbb D}(\chi,\psi;y)$ as measuring 
the distance between the characters $\chi$ and $\psi$ (up to 
some point $y$).  As we will see below (Lemma 3.1) the triangle 
inequality holds:
$$
{\Bbb D}(\chi_1, \psi_1;y) + {\Bbb D}(\chi_2,\psi_2;y) 
\ge {\Bbb D}(\chi_1 \chi_2, \psi_1\psi_2; y). \tag{2.3}
$$
Note that $0\le {\Bbb D}(\chi,\psi;y) \le (1+o(1))\sqrt{2\log \log y}$. 

\proclaim {Definition} Let $\chi$ and $\psi$ be two 
characters and let $\delta >0$.  We say that a character $\chi$ is 
$(\psi,y,\delta)$-pretentious if 
$$
{\Bbb D}(\chi,\psi;y)^2 = \sum_{p\le y} \frac{1-\text{\rm Re }\chi(p) 
\overline{\psi}(p)}{p} \le \delta \log \log y.
$$
\endproclaim

Our main results, from which Theorems 1 and 2 will follow, are 
the following two Theorems.

\proclaim{Theorem 2.1} Of all primitive characters with 
conductor below $(\log q)^{\frac 13}$ let $\xi \pmod m$ 
denote that character for which ${\Bbb D}(\chi,\xi;q)$ 
is a minimum.\footnote{If there are several characters attaining 
this minimum then one can pick $\xi$ to be any one of those characters.}  Then 
$$
M(\chi) \ll (1-\chi(-1)\xi(-1)) \frac{\sqrt{qm}}{\phi(m)} 
\log q \exp\Big( -\frac 12 
{\Bbb D}(\chi,\xi;q)^{2}\Big) + \sqrt{q} (\log q)^{\frac 67}. 
$$
Thus $M(\chi) \ll \sqrt{q}(\log q)^{\frac 67}$ unless 
$\xi(-1)= -\chi(-1)$ and 
$\chi$ is $(\xi,q,\frac 27)$-pretentious.
\endproclaim

In the opposite direction we will show that if $\chi$ is close to 
a character with small conductor (and opposite parity to 
$\chi$) then $M(\chi)$ is large.  

\proclaim{Theorem 2.2}  Let $\psi \pmod \ell$ be a primitive character 
with $\psi(-1)=-\chi(-1)$.  Then 
$$
M(\chi) +\frac{\sqrt{q\ell}}{\phi(\ell)}
 \log \log q \gg \frac{\sqrt{q\ell}}{\phi(\ell)} \log q 
\exp(-{\Bbb D}(\chi,\psi;q)^2).
$$
\endproclaim

We next turn to results conditional on GRH. Given  a real number 
$y \ge 1$ we let ${\Cal S}(y)$ denote the set of integers all of whose 
prime factors are below $y$.  We are motivated by Littlewood's 
conditional result (1.5) which shows that $L(1,\chi)$ is well 
approximated by $\sum_{n\in {\Cal S}(\log^2 q)} \chi(n)/n$.  
We will show that the terms in (2.1) involving 
$e(\pm nx/q)$ may also be replaced by sums involving only
smooth numbers $n$.

\proclaim{Proposition 2.3} Assume GRH.  Let $\chi$ be a primitive 
character $\pmod q$ and let $\alpha$ be a real number. 
Then 
$$
\sum\Sb n\le q \endSb \frac{\cbar(n)}{n} e(n\alpha)= 
\sum\Sb n\le q\\ n\in {\Cal S}((\log q)^{12}) \endSb 
\frac{\cbar(n)}{n} e(n\alpha) + O(1).
$$
\endproclaim 

It follows at once from (1.6), Proposition 2.3 and (2.1) that 
$$
\Big| \sum_{n\le q\alpha}\chi(n)\Big| 
\le \frac{\sqrt{q}}{\pi} |L(1,\cbar)| + \frac{\sqrt{q}}{\pi} \prod_{p\le 
(\log q)^{12}} \Big(1-\frac{1}{p}\Big)^{-1} 
\le (14e^{\gamma}+o(1)) \frac{\sqrt{q}}{\pi}\log \log q.
$$
Thus Proposition 2.3 already gives a refinement of (1.2), and its proof 
(given in \S 5) is simpler than the original proof of (1.2). 
With Proposition 2.3 in place we can argue as in Theorems 2.1 and 
2.2 and arrive at the following conditional analogues from which 
we will deduce Theorems 4 and 5.  

\proclaim{Theorem 2.4} Assume GRH.  Of all primitive characters with 
conductor below $(\log \log q)^{\frac 13}$ let $\xi \pmod m$ 
denote that character for which ${\Bbb D}(\chi,\xi;\log q)$ 
is a minimum.  Then 
$$
M(\chi) \ll (1-\chi(-1)\xi(-1)) \frac{\sqrt{qm}}{\phi({m})} 
\log \log q \exp\Big( -\frac 12 
{\Bbb D}(\chi,\xi;\log q)^{2}\Big) + \sqrt{q} (\log \log q)^{\frac 67}.  
$$
Thus $M(\chi) \ll \sqrt{q}(\log \log q)^{\frac 67}$ 
unless $\xi(-1)= -\chi(-1)$ and 
$\chi$ is $(\xi,\log q,\frac 27)$-pretentious.
\endproclaim

\proclaim{Theorem 2.5} Assume GRH. Let $\psi \pmod \ell$ be a 
primitive character 
with $\psi(-1)=-\chi(-1)$.  Then 
$$
M(\chi) +\frac{\sqrt{q\ell}}{\phi(\ell)}\log \log \log q \gg 
\frac{\sqrt{q\ell}}{\phi({\ell})} 
\log \log q 
\exp(-{\Bbb D}(\chi,\psi;\log q)^2).
$$
\endproclaim

Our work also allows us to make the following refined 
version of Conjecture 1.  

\proclaim{Conjecture 2.6} Let $\chi$ be a primitive character $\pmod q$.  
If $1\leq x\leq q/2$ then $M(\chi) \le (e^{\gamma}/\pi +o(1))\sqrt{q}
\log \log q$ and equality holds here if and only if $
\chi(-1)=-1$, $x \geq q/(\log q)^{o(1)}$,
and 
$$
\sum_{p\le \log q} \frac{1-\text{\rm Re }\chi(p)}{p} = o(1).
$$
Further, if $1\le x\le x+y\le q$, 
$$
\Big| \sum_{x\le n\leq x+y} \chi(n) \Big| \leq \Big(
\frac{2e^\gamma}{\pi \sqrt{3}} +o(1)\Big)   \sqrt{q} \ \log\log
q;
$$
and equality holds here if and only if $\chi(-1)=1$, 
both $|x-q/3|$ and $|x+y-2q/3|$ are $\leq q/(\log q)^{h(q)}$
where $h(q)\to \infty$ as $q\to \infty$, and 
$$
\sum\Sb p\le \log q \\ p\neq 3\endSb \frac{1 
-\text{\rm Re}(\chi(p)\fracwithdelims(){p}{3})}{p}=o(1).
$$
\endproclaim

To finish the paper we use the improved upper bounds for 
$L(1,\chi)$ given in [9], to obtain a modest improvement over 
Hildebrand's results [10] on the constant in 
the P{\' o}lya-Vinogradov inequality. 

\proclaim{Theorem 2.7}  Let $\chi$ be a primitive character 
$\pmod q$, and set $c=1/4$ if $q$ is cubefree, and $c=1/3$ 
otherwise.  If $\chi(-1)=1$ then 
$$
M(\chi) \le \frac{69}{70} 
\frac{c+o(1)}{\pi \sqrt{3}} \sqrt{q} \log q.
$$
If $\chi(-1)=-1$ then 
$$
M(\chi) \le \frac{c+o(1)}{\pi} \sqrt{q} \log q.
$$ 
\endproclaim 

In the case $\chi(-1)=1$ our result improves Hildebrand's estimate
by a factor of $\tfrac {69}{70}$.  Hildebrand gives an estimate 
for a slightly different quantity than $M(\chi)$ when $\chi(-1)=-1$.  

The exponents in Theorems 1, 2, 4, and 5 can all be improved 
by refining the technical Lemmas 3.4 and 4.3. Although we can give 
some improvements to both Lemmas, we have refrained from 
doing so in the interests of a simpler exposition. 
We invite the reader to attain our objectives and reap the 
improved theorems, by replacing the lower bound given in Lemma 3.4 by the 
best possible result, and to replace the exponent `$1/2$' by `$1$' in Lemma 4.3 
(perhaps at the cost of a term of smaller order of magnitude).

\head 3.  The distance between characters and deductions \endhead

\noindent In this section we gather together information on 
the distance between characters defined in (2.2) and show how 
Theorems 1, 2, 4, and 5 may be deduced from Theorems 2.1, 2.2, 2.4 
and 2.5.  Let ${\bold z}$ and ${\bold w}$ denote sequences 
$(z(2),z(3),\ldots)$ and $(w(2),w(3),\ldots)$ indexed by the primes, and 
such that $|z(p)|\le 1$ and $|w(p)|\le 1$ for all $p$.  For two such 
sequences we define (generalizing (2.2)) 
$$
{\Bbb D}({\bold z},{\bold w};y) = \Big(\sum_{p\le y} 
\frac{1-\text{Re }z(p)\overline{w(p)}}{p} 
\Big)^{\frac 12}.
$$
Given two sequences ${\bold z}_1$ and ${\bold z}_2$ we will denote 
by ${\bold z}_1 {\bold z}_2$ the sequence obtained by 
multiplying componentwise: $(z_1(2)z_2(2), z_1(3) z_2(3), \ldots)$.  

\proclaim{Lemma 3.1}  With the above notations we have 
the triangle inequality 
$$
{\Bbb D}({\bold z}_1,{\bold w}_1;y) + {\Bbb D}({\bold z}_2, {\bold w}_2;y) 
\ge {\Bbb D}({\bold z}_1 {\bold z}_2, {\bold w}_1 {\bold w}_2; y).
$$
\endproclaim

\demo{Proof} Since ${\Bbb D}({\bold z},{\bold w};y) = 
{\Bbb D}({\bold 1}, {\bold {\overline{z}w}};y)$ 
we may assume that ${\bold z}_1 ={\bold 1} ={\bold z}_2$.  
Using the Cauchy-Schwartz inequality we see that 
$({\Bbb D}({\bold 1}, {\bold w}_1;y) +{\Bbb D}({\bold 1},{\bold w}_2;y))^2$
is 
$$
\align
&= \sum_{p\le y} \Big(\frac{1-\text{Re }w_1(p)}{p} 
+\frac{1-\text{Re }w_2(p)}{p} 
\Big) + 2 {\Bbb D}({\bold 1},{\bold w}_1;y) 
{\Bbb D} ({\bold 1},{\bold w}_2;y)\\
&\ge 
\sum_{p\le y} \frac 1p 
\Big(1-\text{Re }w_1(p) +1- \text{Re }w_2(p) + 2\sqrt{1-\text{Re }w_1(p)} 
\sqrt{1-\text{Re }w_2(p)} \Big)\\
&\ge \sum_{p\le y} \frac{1}{p} \Big( 1-\text{Re }w_1(p) +1- \text{Re }w_2(p) 
+|\text{Im }w_1(p)||\text{Im } w_2(p)|\Big)\\
&\ge \sum_{p\le y} \frac 1p \Big(1- \text{Re }w_1(p)w_2(p)\Big),\\
\endalign
$$
which proves the Lemma.
\enddemo 

More generally, that given a sequence $(a(2),a(3),\ldots)$ 
of non-negative real numbers we could define the distance between 
${\bold z}$ and ${\bold w}$ as 
$(\sum_{p\le y} a(p) (1-\text{Re }z(p){\overline {w(p)}} ) 
)^{\frac 12}$.  A simple modification of the proof above shows 
that this  also satisfies 
the triangle inequality. 

We now turn to estimates on distances between characters.  
We first record a consequence of the prime number theorem in arithmetic 
progressions that we will find useful below.  Suppose 
$a\pmod \ell$ is a reduced residue class.  Then for any $x$ such 
that $\ell \le (\log x)^A$ ($A$ an arbitrary constant) we have that 
$$
\sum\Sb \ell \le p \le x\\ p \equiv a\pmod \ell \endSb 
\frac 1p = (1+o(1)) \frac{1}{\phi(\ell)} \log \log x. \tag{3.1}
$$  

\proclaim{Lemma 3.2} Let $\chi \pmod q$ be a primitive character 
of odd order $g$.  Suppose $\xi \pmod m$ is a primitive character 
such that $\chi(-1)\xi(-1)=-1$.  If $m\le (\log y)^A$ then 
$$ 
{\Bbb D}(\chi,\xi;y)^2 \ge  (\delta_g +o(1)) \log \log y.
$$
\endproclaim

\demo{Proof} Since $\chi$ has odd order, $\chi(-1)=1$. 
Thus $\xi(-1)=-1$ and $\xi$ must have even order $k\ge 2$ 
say.  We have 
$$ 
{\Bbb D}(\chi,\xi;y)^2  
\ge \sum_{-k/2 <\ell \le k/2} 
\Big( \sum\Sb p\le y\\ \xi(p)=e(\ell/k)\endSb \frac 1p\Big) 
\min_{z^g=0,1} (1-\text{Re } ze(-\ell/k)).
$$
If $\Vert \lambda \Vert$ denotes the distance of $\lambda$ 
from the nearest integer then we may check that 
$\min_{z^g =0, 1} (1-\text{Re }ze(-\ell/k)) 
= 1- \cos(\frac{2\pi}{g} \Vert \ell g/k\Vert)$.  
An application of (3.1) gives that 
$$
\sum\Sb p\le y\\ \xi(p)=e(\ell/k)\endSb \frac 1p \ge 
(1+o(1)) \frac 1k \log \log y.
$$
Writing 
$g/k=g^*/k^*$ in lowest terms (note that $k^* \ge 2$ is 
even) we deduce that 
$$
\align
{\Bbb D}(\chi,\xi;y)^2  &\ge (1+o(1)) \frac 1k \log \log y  \ \cdot \ 
\frac k{k^*} 
\sum_{-k^*/2 <\ell \le k^*/2} (1- \cos\tfrac{2\pi \ell}{gk^*})\\
&\sim \Big( 1- \frac{\sin(\pi/g)}{k^* \tan(\pi/(gk^*))}\Big)\log \log y.
\\
\endalign
$$
Since $k^* \tan(\pi/(gk^*))> \pi/g$, the Lemma follows.
\enddemo 

\demo{Deducing Theorem 1 from Theorem 2.1}  Suppose $\chi$ has 
odd order $g$ and let $\xi$ be the character with 
conductor below $(\log q)^{\frac 13}$ with smallest ${\Bbb D}(\chi,\xi;q)$. 
If $\chi(-1)\xi(-1)=1$ then Theorem 2.1 gives $M(\chi)\ll \sqrt{q} 
(\log q)^{\frac 67}$ which is stronger than 
Theorem 1.  If $\chi(-1)\xi(-1)=-1$ then Lemma 3.2 
gives that ${\Bbb D}(\chi,\xi;q)^2 \ge (\delta_g +o(1)) \log 
\log q$, and Theorem 1 follows at once from Theorem 2.1.

\enddemo 

\demo{Deducing Theorem 4 from Theorem 2.4} This is entirely 
analogous to the above deduction. 

\enddemo

\proclaim{Lemma 3.3} Let $g\ge 2$ be fixed.  
Suppose that for $1\le j\le g$,  
$\chi_j \pmod{q_j}$ is a primitive character.    
Let $y$ be large, and suppose $\xi_j \pmod {m_j}$ are 
primitive characters with conductors $m_j \le \log y$. 
Suppose that $\chi_1\cdots \chi_g$ is the trivial character, 
but $\xi_1\cdots \xi_g$ is not trivial.  Then 
$$
\sum_{j=1}^{g} {\Bbb D}(\chi_j,\xi_j;y)^2 
\ge \Big(\frac{1}{g} +o(1)\Big) \log \log y.
$$
\endproclaim

\demo{Proof} We decompose ${\Bbb D}(\chi_j,\xi_j;y)^2$ as 
${\Bbb D}_0(\chi_j,\xi_j;y)^2 +{\Bbb D}_1(\chi_j,\xi_j;y)^2$ 
where in ${\Bbb D}_0$ we sum over primes $p\le y$ dividing 
the l.c.m. of $q_1$, $\ldots$, $q_g$, and in ${\Bbb D}_1$ 
we sum over all other primes $p\le y$.  Then the triangle inequality 
holds for ${\Bbb D}_1$, and using Cauchy-Schwartz we find 
that 
$$
\sum_{j=1}^{g} {\Bbb D}_1(\chi_j,\xi_j;y)^2 \ge 
\frac 1g \Big(\sum_{j=1}^{g} {\Bbb D}_1(\chi_j,\xi_j;y)\Big)^2 
\ge \frac 1g \sum\Sb p\le y\\ p\nmid q_1\cdots q_g \endSb 
\frac{1-\text{Re }\overline{\xi_1\cdots\xi_g}(p)}{p}.
$$
Trivially 
$$
\sum_{j=1}^{g} {\Bbb D}_0(\chi_j,\xi_j;y)^2 
\ge \sum\Sb p\le y\\ p|q_1\cdots q_g\endSb \frac 1p 
\ge \frac 1g \sum\Sb p\le y\\ p|q_1\cdots q_g\endSb 
\frac{1-\text{Re }\overline{\xi_1\cdots\xi_g}(p)}{p},
$$
and so we deduce that 
$$
\sum_{j=1}^{g} {\Bbb D}(\chi_j,\xi_j;y)^2 \ge 
\frac 1g {\Bbb D}(1,\xi_1\cdots\xi_g;y)^2.
$$
The Lemma now easily follows from (3.1).
\enddemo

\demo{Deducing Theorem 2 from Theorems 2.1 and 2.2} We 
first consider the case when $g$ is odd.  For each $1\le j\le g$ 
let $\xi_j \pmod {m_j}$ denote that primitive 
character with conductor below $(\log q_j)^{\frac 13}$ 
for which ${\Bbb D}(\chi_j,\xi_j;q_j)$ is a minimum.  
If for some $j$ we have $\chi_j(-1)\xi_j(-1)=1$ then 
Theorem 2.1 gives that $M(\chi_j)\ll \sqrt{q_j} (\log q_j)^{\frac 67}$ 
and our claimed bound follows.  Suppose now that 
$\chi_j(-1)\xi_j(-1)=-1$ for all $j$.  By Theorem 2.1, and since
$q_j\le q$, we see 
that 
$$
\align 
\frac{M(\chi_j)}{\sqrt{q_j}} &\ll (\log q_j) \exp\Big(-\frac{1}{2} 
{\Bbb D}(\chi_j,\xi_j;q_j)^2\Big) +(\log q_j)^{\frac 67} 
\\
&\ll (\log q) \exp\Big(-\frac{1}{2} {\Bbb D}(\chi_j,\xi_j;q)^2\Big) 
+ (\log q)^{\frac 67}.\\
\endalign
$$
Therefore 
$$
\prod_{j=1}^{g} \frac{M(\chi_j)}{\sqrt{q_j}} 
\ll (\log q)^{g} \exp\Big(- \frac 12
\sum_{j=1}^{g} {\Bbb D}(\chi_j,\xi_j;q)^2\Big) 
+ (\log q)^{g-\frac 17}.
$$
We know that $\chi_1\cdots \chi_g$ is the 
trivial character, and since $g$ is 
odd $(\xi_1\cdots\xi_g) (-1)=(-1)^g=-1$ 
and so $\xi_1\cdots \xi_g$ is not trivial. 
Lemma 3.1 now gives the bound of the Theorem. 

Now we consider the case when $g$ is even.  
If $g=2$ then $\chi_1$ and $\chi_2$ are 
complex conjugates and there is nothing 
to prove.  Suppose now that $g\ge 4$.  
If for any $1\le j\le g-1$ we have 
$M(\chi_j) \ll \sqrt{q_j} (\log q)^{\frac 67}$ 
then the bound of the Theorem holds trivially.  
Suppose now that for each $1\le j\le g-1$ we have 
that $M(\chi_j)\gg \sqrt{q_j} (\log q)^{\frac 67}$. 
If $\xi_j \pmod{m_j}$ denotes the primitive 
character with conductor below $(\log q_j)^{\frac 13}$ 
with minimum ${\Bbb D}(\chi_j,\xi_j;q_j)^2$ then 
by Theorem 2.1 we have that $\chi_j(-1)\xi_j(-1)=-1$, and 
that 
$$ 
M(\chi_j) \ll \frac{\sqrt{q_j m_j}}{\phi(m_j)} 
(\log q_j) \exp\Big(-\frac 12 {\Bbb D}(\chi_j,\xi_j;q_j)^2\Big) 
\ll \frac{\sqrt{q_jm_j}}{\phi(m_j)} (\log q) \exp\Big(-\frac 12 
{\Bbb D}(\chi_j,\xi_j;q)^2\Big),
$$
so that (for $j\le g-1$)
$$
{\Bbb D}(\chi_j,\xi_j;q)^2 \le 2 \log \Big(\frac{\sqrt{q_jm_j}\log q}
{M(\chi_j)\phi(m_j)}\Big)+O(1).
$$
Now note that $\chi_g$ is the primitive character inducing 
$\overline{\chi_1\cdots \chi_{g-1}}$ and so we let 
$\psi$ denote the primitive character inducing $\overline{\xi_1
\cdots\xi_{g-1}}$.  We note  (using the triangle  
and the Cauchy-Schwartz inequalities, as well as
 an argument as in Lemma 3.3 to 
handle the primes dividing $q_1\cdots q_{g-1} m_1\cdots m_{g-1}$) that
$$
{\Bbb D}(\chi_g,\psi;q)^2 \le (g-1) \sum_{j=1}^{g-1} 
{\Bbb D}(\chi_j,\xi_j;q)^2,
$$
and that $\chi_g(-1)\psi(-1) =(-1)^{g-1} =-1$.  Appealing 
to Theorem 2.2 we obtain the Theorem.

\enddemo

\demo{Deducing Theorem 5 from Theorems 2.4 and 2.5} This 
is entirely analogous to our deduction above.

\enddemo

Finally we record a lemma which will be useful later.

\proclaim{Lemma 3.4} Let $\chi \pmod q$ be a primitive 
character.  Of all primitive characters with conductor below 
$\log y$, suppose that $\psi_j \pmod {m_j}$ ($1\le j\le A$) 
give the smallest distances ${\Bbb D}(\chi,\psi_j;y)$ 
arranged in ascending order.  Then for each $1\le j \le A$ 
we have that 
$$
{\Bbb D}(\chi,\psi_j;y)^2 \ge \Big(1-\frac{1}{\sqrt{j}} +o(1)
\Big) \log \log y.
$$
\endproclaim

\demo{Proof} Notice that 
$$
\align
{\Bbb D}(\chi,\psi_j;y)^2 &\ge \frac{1}{j} \sum_{k=1}^{j} {\Bbb 
D}(\chi,\psi_k;y)^2
= \frac 1j \sum_{p\le y} \frac 1p \sum_{k=1}^{j} (1-\text{Re 
}\chi(p)\overline{\psi_k}(p))\\
&\ge \frac 1j \sum_{p\le y} \frac 1p \Big( j - \Big| \sum_{k=1}^j 
\psi_k(p)\Big|\Big). \tag{3.2}\\
\endalign
$$
By Cauchy-Schwartz we have that 
$$ 
\Big( \sum_{p\le y} \Big| \sum_{k=1}^{j} \psi_k(p)\Big|\Big)^2 
\le \Big(\sum_{p\le y} \frac 1p \Big) \Big(\sum_{p\le y} \frac 1p 
\Big|\sum_{k=1}^{j} 
\psi_k(p)\Big|^2 \Big). \tag{3.3}
$$
The first term in the RHS above is $\sim \log \log y$.  The second term is 
$$
\sum_{p\le y} \frac 1p \Big(j+ \sum\Sb 1\le k , \ell \le j \\ k\neq \ell \endSb 
\frac{\psi_k(p)
\overline{\psi_{\ell}}(p)}{p}\Big) 
\sim j\log \log y,
$$
by appealing to (3.1).  Using these estimates to bound the quantity in (3.3), 
and 
inserting that bound in (3.2) we obtain the Lemma.

\enddemo

\head 4. Preliminary Lemmas \endhead

\noindent Here we collect together some lemmas used below.  
For any character $\chi \pmod q$ we recall the 
{\sl Gauss sum}
$$
\tau(\chi) = \sum_{a\pmod q} \chi(a) e(a/q). \tag{4.1}
$$
It is immediate that if $(b,q)=1$ then 
$$
\sum_{a\pmod q} \chi(a) e(ab/q) = \cbar(b) \tau(\chi). \tag{4.2}
$$ 

\proclaim{Lemma 4.1} Suppose that $\chi \pmod q$ is induced by the 
primitive character $\chi^{\prime} \pmod{q^{\prime}}$.  Then 
$$
\tau(\chi) = \mu(q/q^{\prime}) \chi^{\prime}(q/q^{\prime}) 
\tau(\chi^{\prime}).
$$
If $\chi \pmod q$ is primitive then $|\tau(\chi)|=\sqrt{q}$ 
and (4.2) holds for all integers $b$. 
\endproclaim 
\demo{Proof}  Note that 
$$ 
\tau(\chi) = \sum\Sb a\pmod q \\ (a,q/q^{\prime})=1\endSb 
\chi^{\prime}(a) e(a/q) 
= \sum_{d| (q/q^{\prime})} \mu(d) \chi^{\prime}(d) \sum_{a\pmod {q/d} }
\chi^{\prime}(a) e(ad/q).
$$
The inner sum vanishes unless $d=q/q^{\prime}$ and the first result 
follows.  The second statement is well known, see for example [4].
\enddemo

\proclaim{Lemma 4.2} Let $f$ be a completely 
multiplicative function with $|f(n)|\le 1$ for all $n$. 
Suppose $|\alpha -b/r| \le 1/r^2$ with $(b,r)=1$. 
For any $2\le R\le r$ and any $N\ge Rr$ we have 
$$
\sum_{n\le N} f(n) e(n\alpha) \ll \frac{N}{\log N} + N \frac{(\log R)^{3/2}}
{\sqrt{R}}, 
$$
and 
$$ 
\sum_{Rr\le n\le N} \frac{f(n)}{n} e(n\alpha) 
\ll \log \log N + \frac{(\log R)^{3/2}}{\sqrt{R}} \log N. 
$$
\endproclaim
\demo{Proof}  The first bound follows from Corollary 1 of Montgomery 
and Vaughan [12].  The second estimate follows easily from the 
first and partial summation. 

\enddemo 

\proclaim{Lemma 4.3}  If $f$ is a multiplicative 
function with $|f(n)|\le 1$ for all $n$ then 
$$
\sum_{n\le x} \frac{f(n)}{n} \ll 1+ \log x \exp\Big(-\sum_{p\le x} 
\frac{2-|1+f(p)|}{p}\Big)
\ll 1 +\log x \exp\Big(-\frac 12 {\Bbb D}(1,f;x)^2\Big).
$$
Further if $y\ge 1$ then 
$$
\sum\Sb n\le x\\ n\in{\Cal S}(y)\endSb 
\frac{f(n)}{n} \ll 1+ \log y 
\exp\Big(-\frac 12 {\Bbb D}(1,f;y)^2\Big). 
$$
\endproclaim
\demo{Proof} For the first assertion see the remark after Proposition 8.1 of [7] 
and 
note that if $|z|\le 1$ then $2-|1+z| \ge \tfrac 12 (1-\text{Re }z)$. 
To see the second assertion note that if $y\le x$ then 
$$
\align
\sum\Sb n\le x\\ n\in {\Cal S}(y)\endSb 
\frac{f(n)}{n} 
&\ll 1+ \log x \exp\Big(-\sum_{p\le y}\frac{2-|1+f(p)|}{p} -\sum_{y<p\le x} 
\frac 1p\Big) 
\\
&\ll 1+ \log y \exp\Big(-\frac 12 {\Bbb D}(1,f;y)^2\Big). \\
\endalign
$$
If $y>x$ then $\log x \exp(-\frac 12 
{\Bbb D}(1,f;x)^2) \ll \log y \exp(-\frac 12 {\Bbb D}(1,f;y)^2)$, and 
so the second assertion holds in this case also. 
\enddemo 

\proclaim{Lemma 4.4} Let $f$ be a completely multiplicative 
function with $|f(n)|\le 1$ for all $n$.  
Then for any integer $\ell \ge 1$ we have 
$$
\sum_{n\le x} \frac{f(n)}{n} 
= \prod_{p|\ell} \Big( 1-\frac{f(p)}{p}\Big)^{-1} 
\sum\Sb n\le x\\ (n,\ell)=1\endSb \frac{f(n)}{n} + O((\log \log (\ell+2))^2).
$$
\endproclaim
\demo{Proof} Writing $n$ as $uv$ where $u$ is composed only 
of primes dividing $\ell$ and $v$ is coprime to $\ell$ we 
see that 
$$
\sum_{n\le x} \frac{f(n)}{n} = \sum\Sb u \\ p|u\implies 
p|\ell\endSb \frac{f(u)}{u}\sum\Sb v\le x/u\\ (v,\ell)=1\endSb \frac{f(v)}{v}
= \prod_{p|\ell}\Big(1-\frac{f(p)}{p}\Big)^{-1} 
\sum\Sb v\le x\\ (v,\ell)=1\endSb \frac{f(v)}{v} 
+ O\Big(\sum\Sb u \\ p|u\implies 
p|\ell\endSb \frac{\log u}{u}\Big).
$$
The error term is seen to be $\ll (\sum_{p|\ell} \frac{\log p}{p} ) 
\prod_{p|\ell} (1+1/p) \ll (\log \log (\ell +2))^2$.  

\enddemo

\head 5.  Proof of Proposition 2.3 \endhead

\noindent We begin by recalling a consequence of GRH.  If $\psi \pmod m$ 
is a non-principal character then for all $x\ge 2$ 
$$
\sum_{n\le x} \psi(n)\Lambda(n) \ll \sqrt{x} \log x \log (mx). 
$$
This follows from standard arguments: for example, take $T=x^2$ in 
(13) on page 120 of [4] and use GRH.  
It follows from the above and partial summation that 
$$
\sum_{p\le x} \psi(p) \ll \sqrt{x} \log (mx). \tag{5.1}
$$

\proclaim{Lemma 5.1}  Assume GRH.  If $\chi$ is a 
primitive character $\pmod q$ and $x < q^{3/2}$ then, 
uniformly for all $\theta$, we have that 
$$
\sum_{p\le x} \chi(p) e(p\theta) \ll x^{5/6} \log q.
$$
\endproclaim 
\demo{Proof}  First we show that if $(b,r)=1$ with 
$r<q$ then for any $x\ge 2$ 
$$
\sum_{p\le x} \chi(p) e(bp/r) \ll \sqrt{rx} \log (qx). \tag{5.2} 
$$
To see this, note that 
$$
\align
\sum_{p\le x} \chi(p)e(bp/r) 
&= \sum\Sb p\le x\\ (p,r)=1\endSb \chi(p) e(bp/r)
+ O\Big(\sum_{p|r} 1\Big) \\
&= \frac{1}{\phi(r)} \sum_{\psi \pmod r} \overline{\psi}(b) \tau(\psi) 
\sum_{p\le x} \chi(p)\overline{\psi}(p) + O(\log q).\\ 
\endalign
$$
Since $r<q$ and $\chi$ is primitive we know that $\chi\overline{\psi}$ is a 
non-principal character $\pmod {qr}$.  Appealing now to (5.1) 
and using that $|\tau(\psi)|\le \sqrt{r}$ from Lemma 4.1 we obtain (5.2). 

We now turn to the proof of the lemma.  
Set $R=x^{2/3}$ and find $r\le R$ such that 
$\theta=b/r+\beta$ where $(b,r)=1$ and $|\beta|\le 1/(rR)$.  
If $x <q^{3/2}$ then $r<q$ and by (5.2) we obtain that 
$\sum_{p\le N} \chi(p) e(bp/r) \ll \sqrt{rN} \log (qN)$ for all $N\ge 2$. 
By partial summation we see that 
$$
\sum_{p\le x} \chi(p)e(p\theta) = 
\sum_{p\le x} \chi(p)e(pb/r) e(p\beta) 
= \int_2^x e(t\beta) d\Big( \sum_{p\le t} \chi(p) e(pb/r)\Big),
$$
and integrating by parts using our bound above, we obtain that 
$$
\sum_{p\le x} \chi(p) e(p\theta) \ll (1+|\beta|x) \sqrt{rx} \log q 
\ll x^{5/6} \log q.
$$

\enddemo

\proclaim{Lemma 5.2} Assume GRH.  If $\chi$ is a 
primitive character $\pmod q$ and $x<q^{3/2}$ then 
$$ 
\sum_{n\le x} \chi(n) e(n\alpha) 
= \sum\Sb n\le x \\ n\in {\Cal S}(y)\endSb 
\chi(n)e(n\alpha) + O(x y^{-1/6} \log q). 
$$
\endproclaim  

\demo{Proof}  Write $n\notin {\Cal S}(y)$ as $pm$ where $p$ is 
the largest prime divisor of $n$.  Thus $x/m \ge p>y$ and $m \le x/y$ and 
so (denoting by $P(m)$ the largest prime factor of $m$)
$$
\sum\Sb n\le x\\ n\notin {\Cal S}(y) \endSb 
\chi(n) e(n\alpha) 
= \sum\Sb m\le x/y\endSb \chi(m) \sum_{\max(P(m)-1,y) < p\le x/m}  
\chi(p) e(pm\alpha),
$$
and by Lemma 5.1 this is 
$$
\ll \sum_{m\le x/y} (x/m)^{5/6} \log q \ll x y^{-1/6} \log q,
$$
as required.
\enddemo

Lemma 5.2 and partial summation give that 
$$
\sum_{n\le q} \frac{\cbar(n)e(n\theta)}{n} = 
\sum\Sb n\le q\\ n\in {\Cal S}(y)\endSb \frac{\cbar(n)e(n\theta)}{n} 
+O(y^{-1/6} \log^2 q),
$$
and Proposition 2.3 follows.  

\medskip

We remark that Lemma 5.2 with $\alpha=0$ shows how character 
sums may be approximated by character sums involving only 
smooth numbers.  This question is explored in greater 
depth in our paper [6].

\head 6.  Proof of Theorems 2.1, 2.2, 2.4 and 2.5 \endhead

\noindent The main ideas of our proof work whether or not GRH is 
assumed, the only difference being that the relevant parameters 
need to be chosen differently in each case.  To present this in 
a unified manner we adopt the following convention.  We set 
$Q= \log q$ if GRH is assumed, and $Q=q$ if no assumption is 
being made.  Accordingly we warn the reader that the results in 
this section must all be read keeping this convention in mind.

By (2.1), to understand $M(\chi)$ we must gain an 
understanding of $\sum_{n\le q} \cbar(n)e(n\alpha)/n$ 
where $\alpha \in [0,1]$.  If we assume GRH then Proposition 2.3 
shows that we may restrict ourselves to 
$\sum_{n\le q, n\in {\Cal S}((\log q)^{12})} \cbar(n) e(n\alpha)/n$.  
Thus, with our convention, we seek to understand 
$$
\sum\Sb n\le q \\ n\in {\Cal S}(Q^{12}) \endSb \frac{\cbar(n)}{n} 
e(n\alpha), \tag{6.1}
$$
since in the unconditional case the criterion $n\in {\Cal S}(Q^{12})$ is 
vacuous.

We now define $s=(\log Q)^{\frac 13}$ and $S= \exp((\log Q)^{\frac 56})$. 
We say that $\alpha$ lies on a {\sl minor arc} if there is a rational 
approximation $|\alpha -b/r|\le 1/(rS)$ with $(b,r)=1$ and $s<r\le S$.  
Otherwise 
we say that $\alpha$ lies on a {\sl major arc}; in this case there 
is a rational approximation $|\alpha-b/r|\le 1/(rS)$ with 
$(b,r)=1$ and $r\le s$.

\proclaim{Lemma 6.1} With the above conventions, if $\alpha$ lies 
on a minor arc then 
$$
\sum\Sb n\le q\\ n\in {\Cal S}(Q^{12}) \endSb 
\frac{\cbar(n)}{n} e(n\alpha) 
\ll (\log Q)^{\frac 56 +o(1)}.
$$
\endproclaim
\demo{Proof}  Suppose $|\alpha-b/r|\le 1/(rS)$ 
where $(b,r)=1$ and $s<r\le S$.  

First we consider the unconditional 
case.  By Lemma 4.2 with $R=r$ we  see that 
$$
\sum_{n\le q}\frac{\cbar(n)e(n\alpha)}{n} 
= \sum_{n\le r^2} \frac{\cbar(n)e(n\alpha)}{n}
+ \sum_{r^2 \le n\le q} \frac{\cbar(n)e(n\alpha)}{n} 
\ll (\log q)^{\frac 56 +o(1)},  
$$
which proves the Lemma in this case.

Now we consider the case when we assume GRH.  
By Lemma 4.2 with $R=r$  we see that
$$
\sum\Sb r^2 \le n\le (\log q)^{\log s}\\ 
n\in {\Cal S}((\log q)^{12}) \endSb
\frac{\cbar(n)}{n}e(n\alpha) 
\ll \log \log s + \log\log\log q
+ \frac{(\log s)^{5/2}}{\sqrt{s}} \log \log q.
$$
Further 
$$
\Big|\sum\Sb n> (\log q)^{\log s}\\ n\in{\Cal S}((\log q)^{12}) \endSb 
\frac{\cbar(n)}{n} e(n\alpha)\Big| 
\le \frac{1}{s} \sum_{n\in {\Cal S}((\log q)^{12})}
 \frac{1}{n^{1-1/\log \log q}} 
\ll \frac{1}{s} \log \log q,
$$
and, trivially, 
$$
\sum\Sb n\le r^2 \\ n\in {\Cal S}((\log q)^{12})\endSb 
\frac{\cbar(n)}{n} e(n\alpha) \ll \log r \le \log S.
$$
Combining these estimates we get the Lemma in this 
situation.
\enddemo

We now consider (6.1) when $\alpha$ lies on a major arc.  
Thus we suppose that $|\alpha-b/r| \le 1/(rS)$ where $(b,r)=1$ and $r\le s$, 
and no such approximation exists with $s\le r\le S$. Define 
$N=N_{q,\alpha,b/r} = \min(q,1/|r\alpha-b|)$.

\proclaim{Lemma 6.2}   With the above conventions, we have 
$$
\sum\Sb n\le q\\ n\in {\Cal S}(Q^{12}) \endSb 
\frac{\cbar(n)}{n} e(n\alpha) = 
\sum\Sb n\le N \\ n\in {\Cal S}(Q^{12}) \endSb 
\frac{\cbar(n)}{n} e(nb/r) + O(\log \log Q). 
$$
\endproclaim 

\demo{Proof}  If $N=q$ then the lemma follows easily from 
$|e(n\alpha)-e(nb/r)| \ll n|\alpha-b/r| \le n/N$. Now 
suppose that $S\le N=1/|r\alpha-b| <q$.  We find an approximation 
$|\alpha-b_1/r_1| \le 1/(r_1N)$ where $(b_1,r_1)=1$ and 
$r_1\le N$.  Note that 
$1/(r r_1) \le |b/r-b_1/r_1| \le 1/(rN)+1/(r_1N)$ and 
so $r_1 \ge N-r \ge N-s \ge N/2$.  We now set $R=(\log Q)^5$ 
and divide the interval $(N,q]$ into three intervals:  
$I_1$ which contains the integers in $(N,q]$ that 
are in $(N,Rr_1]$, $I_2$ which contains the integers in $(N,q]$ 
that are in $(Rr_1,\exp((\log Q)^2)]$, and $I_3$ which 
contains the integers in $(N,q]$ that are larger than 
$\exp((\log Q)^2)$.  

Since $N/2 \le r_1 \le N$ it follows that  
$$ 
\sum\Sb n\in I_1 \\ n\in {\Cal S}(Q^{12})\endSb 
\frac{\cbar(n)}{n} e(n\alpha) \ll \log R \ll \log \log Q.   
$$
An application of Lemma 4.2 shows that 
$$
\sum\Sb n\in I_2 \\ n\in {\Cal S}(Q^{12})\endSb 
 \frac{\cbar(n)}{n} e(n\alpha) \ll \log \log Q. 
$$
Finally, since each element of $I_3$ is at least 
$\exp((\log Q)^2)$ we see that 
$$ 
\sum\Sb n\in I_3 \\ n\in {\Cal S}(Q^{12}) \endSb 
\frac{\cbar(n)}{n} e(n\alpha) 
\ll \frac{1}{Q} \sum_{n\in {\Cal S}(Q^{12})} \frac{1}{n^{1-1/\log Q}} 
\ll 1.
$$
Combining these estimates we obtain that 
$$
\sum\Sb n\le q\\ n\in {\Cal S}(Q^{12}) \endSb 
 \frac{\cbar(n)}{n} e(n\alpha) 
= \sum\Sb n\le N \\n\in{\Cal S}(Q^{12})\endSb
 \frac{\cbar(n)}{n} e(n\alpha) + O(\log \log Q),
$$
and since $|e(n\alpha)- e(nb/r)| \ll n|\alpha-b/r| \le n/N$, the Lemma 
follows.

\enddemo

\subhead 6.1. Lower bounds for $M(\chi)$: 
Proof of Theorems 2.2 and 2.5 \endsubhead
 
\noindent We consider the quantity (6.1) 
for $\alpha_{b,N}= b/\ell +1/N$ where $b$ runs over 
reduced residue classes $\pmod \ell$ and 
$1\le N\le q$.  We multiply this by $\overline{\psi}(b)$ 
and sum over all reduced residue classes $b \pmod \ell$. 
Thus we arrive at 
$$
\sum_{b\pmod \ell} \overline{\psi}(b)
\sum\Sb n\le q\\ n\in {\Cal S}(Q^{12}) \endSb 
\frac{\cbar(n)}{n} e(n\alpha_{b,N}) 
= \tau(\overline{\psi}) \sum\Sb n\le q\\ n\in {\Cal S}(Q^{12}) 
\endSb \frac{(\cbar\psi)(n)}{n} e(n/N).
$$
Exactly as in the proof of Lemma 6.2, set $R=(\log Q)^{5}$ 
and divide the integers in $(N,q]$ into 
intervals $I_1$, $I_2$ and $I_3$.  Then we deduce 
that 
$$
\sum_{b\pmod \ell} \overline{\psi}(b)
\sum\Sb n\le q\\ n\in {\Cal S}(Q^{12}) \endSb 
\frac{\cbar(n)}{n} e(n\alpha_{b,N}) 
= \tau(\overline{\psi}) \sum\Sb n\le N \\ 
n\in {\Cal S}(Q^{12}) \endSb \frac{(\cbar\psi)(n)}{n} 
+ O(\sqrt{\ell} \log \log Q).
$$
Now consider $\sum_{b=1}^{\ell} \overline{\psi}(b) 
\sum_{n\le q\alpha_{b,N}} \chi(n)$ which in magnitude 
is plainly $\le \phi(\ell) M(\chi)$.  We see by (2.1), 
Proposition 2.3 (in the conditional case), and the above 
remarks that if $\ell >1$
$$
\sum_{b=1}^{\ell} \overline{\psi}(b) 
\sum_{n\le q\alpha_{b,N}} \chi(n) 
= -\frac{\tau(\chi)\tau(\overline{\psi})}{2\pi i}(\psi(-1)-\chi(-1))
\sum\Sb n\le N \\ 
n\in {\Cal S}(Q^{12}) \endSb \frac{(\cbar\psi)(n)}{n} 
+O(\sqrt{q\ell }\log \log Q), \tag{6.2a}
$$
while if $\ell=1$ we have (because of the extra $L(1,\cbar)$ term)
$$
\sum_{n\le q\alpha_{1,N}} \chi(n) 
= \frac{\tau(\chi)}{2\pi i} (1-\chi(-1)) 
\sum\Sb N\le n\le q \\ n \in {\Cal S}(Q^{12}) \endSb 
\frac{\cbar(n)}{n} +O(\sqrt{q} \log \log Q). \tag{6.2b}
$$

If we assume GRH then taking $N=q$ in the 
case $\ell >1$ and $N=1$ when $\ell =1$ 
we obtain easily that the sum over 
$n$ in (6.2a,b) is $\gg \log Q \exp(-{\Bbb D}(\chi,\psi;Q)^2)$
and Theorem 2.5 follows.  In the unconditional 
case we show in the next Lemma that a similar lower 
bound holds for some $1\le N \le q$ which proves 
Theorem 2.2.

\proclaim{Lemma 6.3}  Let $\eta \pmod r$ be a primitive 
character.  Then there exists $1\le N\le r$ such that 
$$
\Big| \sum_{n\le N} \frac{\eta(n)}{n} \Big| +1 \gg \log r 
\exp(-{\Bbb D}(\eta,1;r)^2).
$$
There also exists $1\le N\le r$ with 
$$
\Big| \sum_{N\le n\le r} \frac{\eta(n)}{n} \Big| +1 \gg \log r 
\exp(-{\Bbb D}(\eta,1;r)^2).
$$
\endproclaim

\demo{Proof} Set $\delta= 1/\log r$ and 
observe that 
$$
\frac{1}{e} \sum_{n\le r} \frac{\eta(n)}{n} + 
\int_1^r \frac{\delta}{t^{1+\delta}} 
\sum_{n\le t} \frac{\eta(n)}{n} dt =\sum_{n\le r} 
\frac{\eta(n)}{n^{1+\delta}}.
$$
It follows that 
$$
\max_{N\le r} \Big|\sum_{n\le N} \frac{\eta(n)}{n} \Big| 
\ge \Big| \sum_{n\le r} \frac{\eta(n)}{n^{1+\delta}}\Big|.
$$
We see easily that $L(1+\delta,\eta)= \sum_{n\le r} \eta(n)/n^{1+\delta} 
+O(1)$, and from the Euler product 
that $L(1+\delta,\eta)\gg \log r  \exp(-{\Bbb D}(\eta,1;r)^2)$.  
The first part of the Lemma follows.  The second part is 
similar starting from 
$$
 \sum_{n\le r} \frac{\eta(n)}{n} 
-\int_1^r \frac{\delta}{t^{1+\delta}} 
\sum_{t\le n\le r} \frac{\eta(n)}{n} dt =\sum_{n\le r} 
\frac{\eta(n)}{n^{1+\delta}}.
$$

\enddemo

\subhead 6.2. Upper bounds for $M(\chi)$: Proof 
of Theorems 2.1 and 2.4 \endsubhead

\noindent We continue from Lemma 6.2 our analysis of (6.1) in 
the case when $\alpha$ lies on a major arc.  Of all 
characters with conductor below $s$ we let $\xi \pmod m$ 
denote that character for which ${\Bbb D}(\chi,\psi;Q)$ 
is a minimum.

\proclaim{Lemma 6.4} We keep the conventions of this 
section.  Suppose $(b,r)=1$ with $r\le s$.  
Then 
$$
\sum\Sb n\le N \\ n\in {\Cal S}(Q^{12})\endSb 
\frac{\cbar(n)}{n}e(nb/r) = O((\log Q)^{\frac 67}),
$$
unless $m|r$ in which case it equals 
$$
\frac{\xi(b) \tau(\overline{\xi}) }{\phi(r)} 
\prod\Sb p^{\alpha} \parallel r/m\\ \alpha \ge 1\endSb 
 (\cbar(p^{\alpha}) 
-\overline{\xi}(p) \cbar(p^{\alpha-1})) 
\sum\Sb n\le N\\ n\in {\Cal S}(Q^{2}) \endSb 
\frac{(\cbar\xi)(n)}{n}+ O( (\log Q)^{\frac 67}) .
$$ 
\endproclaim 

\demo{Proof} Note that 
$$
\sum\Sb n\le N\\ n\in {\Cal S}(Q^{12}) \endSb 
\frac{\cbar(n)}{n} e(nb/r) = 
\sum_{d|r} \frac{\cbar(d)}{d} \sum\Sb n\le N/d\\ n \in {\Cal S}(Q^{12}) 
\\ (n,r/d)=1\endSb \frac{\cbar(n)}{n} e\Big(\frac{nb}{r/d}\Big). \tag{6.3}
$$
Since $(nb,r/d)=1$ we see that 
$$
\align
e\Big(\frac{nb}{r/d}\Big) &= \frac{1}{\phi(r/d)} \sum_{a \pmod {r/d}} 
e\Big(\frac{a}{r/d}\Big) \sum_{\psi \pmod{r/d}} \overline{\psi}(a){\psi}(nb)\\ 
&= \frac{1}{\phi(r/d)} \sum_{\psi \pmod{r/d}}{\psi}(nb) \tau(\overline{\psi}).
\\
\endalign
$$
Therefore 
$$
 \sum\Sb n\le N/d\\ n \in {\Cal S}(Q^{12}) 
\\ (n,r/d)=1\endSb \frac{\cbar(n)}{n} e\Big(\frac{nb}{r/d}\Big)= 
\frac{1}{\phi(r/d)} \sum_{\psi \pmod{r/d}} \tau(\overline{\psi}) 
{\psi}(b) 
\sum\Sb n\le N/d \\ n\in {\Cal S}(Q^{12})\endSb 
\frac{(\cbar \psi)(n)}{n}.  \tag{6.4} 
$$
 
By Lemma 4.3 we see that
$$
\Big|\sum\Sb n\le N/d \\ n\in {\Cal S}(Q^{12}) \endSb 
\frac{(\cbar \psi)(n)}{n} \Big| \ll 1 + (\log Q) \exp\Big( -\frac{1}{2} {\Bbb 
D}(\chi,\psi;Q)^2
\Big). 
$$
Using Lemma 3.4 we see that if $\psi$ is not induced by $\xi$ 
then 
$$
{\Bbb D}(\chi,\psi;Q)^2 \ge (1-1/\sqrt{2}+o(1)) \log \log Q,
$$ 
and further that there are at most $9$ characters $\psi \pmod {r/d}$ 
for which ${\Bbb D}(\chi,\psi;Q)^2 \le \frac 23 \log \log Q$.  
Since $|\tau(\overline{\psi})|\le \sqrt{r/d}$ we 
deduce that the contribution of all characters not induced by $\xi$ to 
(6.4) is 
$$
\ll \frac{\sqrt{r/d}}{\phi(r/d)} 
(\log Q)^{ \frac 12 +\frac{1}{2\sqrt{2}} +o(1)} + 
\sqrt{r/d} (\log Q)^{\frac 23}.
$$
The contribution of these terms to (6.3) is 
$\ll (\log Q)^{\frac 12+\frac{1}{2\sqrt{2}}+o(1)} 
+ \sqrt{r} (\log Q)^{\frac 23+o(1)} 
\ll (\log Q)^{\frac 67}$.

We must now handle the contribution to (6.3) from 
characters $\psi$ induced by $\xi \pmod m$.  
If $m\nmid r$ then there 
are no such characters $\psi$, and the Lemma follows in 
this case.  If $m|r$ then we have to account for the 
contribution of the characters $\psi \pmod {r/d}$ induced by $\xi$ 
(thus $d$ must be a divisor of $r/m$).  By Lemma 4.1 and (6.4) we see that 
the contribution of these induced characters to (6.3) is
$$
\sum_{d|r/m} \frac{\cbar(d)}{d} \frac{1}{\phi(r/d)} 
\xi(b) \tau(\overline{\xi}) \mu\Big(\frac{r}{dm}\Big) 
\overline{\xi}\Big(\frac{r}{dm}\Big) \sum\Sb n\le N/d \\ 
(n,r/d)=1\\ n\in {\Cal S}(Q^{12})\endSb 
\frac{(\cbar \xi)(n)}{n}. \tag{6.5}
$$

By Lemma 4.4 
$$ 
\align 
\sum\Sb n\le N/d \\ (n,r/d)=1\\ n\in {\Cal S}(Q^{12}) \endSb
\frac{(\cbar\xi)(n)}{n} 
&=\sum\Sb n\le N \\ (n,r/d)=1\\ n\in {\Cal S}(Q^{12}) \endSb
\frac{(\cbar\xi)(n)}{n} +O(\log d) \\ 
&= \prod_{p|r/d} \Big( 1-\frac{(\cbar \xi)(p)}{p} 
\Big) \sum\Sb n\le N\\ n\in {\Cal S}(Q^{12})\endSb 
 \frac{(\cbar\xi)(n)}{n} + O(\log \log Q).
\\
\endalign
$$
Therefore (6.5) equals, up to an 
error $O(\log \log Q)$, 
$$
\xi(b) \tau(\overline{\xi}) 
\sum_{d|r/m} 
\frac{\cbar(d)}{d} \frac{1}{\phi(r/d)} 
\mu\Big(\frac{r}{md}\Big) \overline{\xi}\Big(\frac{r}{md}\Big) 
\prod_{p| r/(md)} \Big(1-\frac{(\cbar\xi)(p)}{p}\Big) 
\sum\Sb n\le N\\ n\in {\Cal S}(Q^{12}) \endSb 
\frac{(\cbar\xi)(n)}{n},  
$$
which by a straight-forward calculation is 
$$
\frac{\xi(b) \tau(\overline{\xi}) }{\phi(r)} 
\prod\Sb p^{\alpha} \parallel r/m\\ \alpha \ge 1\endSb 
 (\cbar(p^{\alpha}) 
-\overline{\xi}(p) \cbar(p^{\alpha-1})) 
\sum\Sb n\le N\\ n\in {\Cal S}(Q^{12}) \endSb 
\frac{(\cbar\xi)(n)}{n}. \tag{6.6} 
$$

To complete the proof of the Lemma, it remains 
to show that the terms in the 
sum in (6.6) may be restricted to $n\in {\Cal S}(Q^{2})$ 
with an acceptable error, in the case where GRH is assumed.  
We must therefore estimate the contribution of terms $n$ which lie in 
${\Cal S}(Q^{12})$ 
but not in ${\Cal S}(Q^2)$.  We may write 
such $n$ uniquely as $p\ell$ where $p$, the largest prime factor 
of $n$, lies between $\max(P(\ell)-1,Q^2)$ and 
$\min(Q^{12},N/\ell)$, 
and $\ell \le N/Q^{2}$ is in ${\Cal S}(Q^{12})$.  
Thus the contribution of these terms is 
$$
\sum\Sb \ell \le N/Q^2 \\ \ell 
\in {\Cal S}(Q^{12})\endSb 
\frac{(\cbar\xi)(\ell)}{\ell} \sum\Sb \max(P(\ell)-1,Q^2) \le p\\ 
p\le \min(Q^{12},N/\ell)\endSb \frac{(\cbar \xi)(p)}{p}.
\tag{6.7} 
$$
Using (5.1) and partial summation to handle 
the primes larger than $Q^2 (\log Q)^2$, and 
estimating the smaller primes trivially, we obtain 
that the sum over primes in (6.7) above is 
$\ll (\log \log Q)/\log Q$. 
Thus  the quantity (6.7) is $\ll \log \log Q$, and the 
Lemma follows.  
\enddemo 

Combining Lemmas 6.1, 6.2 and 6.4 we arrive at the following:

\proclaim{Lemma 6.5} Keep the conventions of this section.  
Then 
$$
\sum\Sb n\le q\\ 
n\in {\Cal S}(Q^{12})\endSb 
\frac{\cbar(n)}{n} 
(e(-n\alpha)-\chi(-1)e(n\alpha))  \ll (\log Q)^{\frac 67}
$$
unless $\alpha$ lies on a major arc $|\alpha -b/r| \le 1/(rS)$ with $r\le s$, 
$(b,r)=1$ and $m|r$, in which case it equals, up 
to an error $O((\log Q)^{\frac 67})$,
$$
\frac{(\xi(-1)-\chi(-1))\xi(b) \tau(\overline{\xi})}{\phi(r)} 
\prod\Sb p^a \parallel r/m \\ a\ge 1\endSb 
(\cbar(p^a)-\overline{\xi}(p) \cbar(p^{a-1})) 
\sum\Sb n\le N \\ n\in {\Cal S}(Q^2) 
\endSb 
\frac{(\cbar \xi)(n)}{n},
$$
where $N =\min (q,1/|r\alpha-b|)$.  
\endproclaim

\demo{Proof of Theorems 2.1 and 2.4} From Lemma 6.5 we 
see easily that 
$$
\align
\sum\Sb n\le q\\ 
n\in {\Cal S}(Q^{12})\endSb 
\frac{\cbar(n)}{n} 
(e(-n\alpha)&-\chi(-1)e(n\alpha))  \ll (\log Q)^{\frac 67} \\
&+ \frac{(1-\chi(-1)\xi(-1))\sqrt{m}}{\phi(m)} 
\max_{N\le q} \Big|\sum\Sb n\le N\\ n\in{\Cal S}(Q^2)\endSb 
\frac{(\cbar\xi)(n)}{n}\Big|. \\
\endalign
$$
Using Lemma 4.3 this is 
$$
\ll (\log Q)^{\frac 67} + 
\frac{(1-\chi(-1)\xi(-1))\sqrt{m}}{\phi(m)} (\log Q) \exp\Big(-\frac 12 
{\Bbb D}(\chi,\xi;Q)^2\Big). \tag{6.8}
$$
We also note that (using $L(1,\cbar)=\sum_{n\le q}\cbar(n)/n+O(1)$ 
and Lemma 4.3 in the unconditional 
case, and Littlewood's (1.5) in the conditional case)
$$
L(1,\cbar) \ll (\log Q) \exp\Big( -\frac 12 {\Bbb D}(\chi,1;Q)^2\Big). 
\tag{6.9}
$$

Using (6.8) and (6.9) in (2.1) (and using Proposition 2.3 in 
the conditional case) we immediately 
obtain Theorems 2.1 and 2.4 in the case when 
$m=1$ and $\xi$ is the trivial character $\xi(n)=1$ 
for all $n$.  In the case when $m>1$ it follows 
from Lemma 3.4 that ${\Bbb D}(\chi,1;Q)^2\ge (1-1/\sqrt{2}+o(1)) 
\log \log Q$.  Using this in (6.9) we obtain the 
bounds claimed in  Theorems 2.1 and 2.4.

\enddemo

\head 7.  Proof of Theorem 6 \endhead
 
\noindent To prove Theorem 6, we assume GRH 
and continue the analysis of the previous 
section (note that $Q=\log q$).  We distinguish 
two cases: when the nearest character $\xi \pmod m$ 
is the trivial character ($m=1$ and $\xi(n)=1$ 
for all $n$), and when $m>1$.  

We start with the easier second case.  By Lemma 
3.4 we have that ${\Bbb D}(\chi,1;\log q)^2 \ge 
(1-1/\sqrt{2}+o(1)) \log \log \log q$ 
and therefore, by (6.9), $L(1,\cbar)= o(\log \log q)$. 
>From this, (2.1), Proposition 2.3, and Lemma 6.5 we obtain that 
$M(\chi)+o(\sqrt{q}\log \log q)$ is 
$$
\le \frac{\sqrt{q}|\tau({\overline{\xi}})|}{\pi} 
\max_{m|r, N\le q} \Big| 
\frac{1}{\phi(r)} \prod\Sb p^{a}\parallel r/m\\ a\ge 
1\endSb (\cbar(p^a) -\overline{\xi}(p) \cbar(p^{a-1}) 
) \sum\Sb n\le N \\ n\in {\Cal S}(\log^2 q) 
\endSb \frac{(\cbar \xi)(n)}{n}\Big|.
$$
By Lemma 4.4 we see that, up to $o(\sqrt{q}\log \log q)$, 
the above is 
$$
\le \frac{\sqrt{qm}}{\pi} 
\max_{m|r, N\le q} \Big| 
\frac{1}{\phi(r)} \prod\Sb p^{a}\parallel r/m\\ a\ge 
1\endSb
\frac{\cbar(p^a)-\overline{\xi}(p) \cbar(p^{a-1})}{1-(\cbar \xi)(p)/p}
\sum\Sb v\le N\\ v\in {\Cal S}(\log^2 q) \\ 
(v,r)=1\endSb  \frac{(\cbar \xi)(v)}{v}\Big|. 
$$
The product above is 
bounded in magnitude by $\prod_{p | r/m} 2p/(p+1)$,  
and the sum over $v$ above has 
size $\le \prod_{p\le \log^2 q, p\nmid r} (1-1/p)^{-1} 
= (2e^{\gamma}+o(1)) (\phi(r)/r)\log \log q$. 
It follows readily that when $m>1$ 
$$
M(\chi) \le \Big(\frac{2e^{\gamma}}{\pi\sqrt{m}} 
+o(1)\Big) \sqrt{q}\log \log q.
$$
Since there is no primitive character $\pmod 2$ we 
have that $m\ge 3$ and so the bounds of Theorem 6 
follow. 

We now consider the more involved case when $m=1$ and 
$\xi$ is the trivial character.  
We consider $\sum_{\alpha_1 q\le n\le \alpha_2 q} 
\chi(n)$ where $0\le \alpha_1 <\alpha_2 \le 1$ and 
by (2.1) and Proposition 2.3 this is 
$$ 
-\frac{\tau(\chi)}{2\pi i}  \sum\Sb n\le q \\ 
n\in {\Cal S}((\log q)^{12}) \endSb 
\frac{\cbar(n)}{n} (e(-n\alpha_1) -\chi(-1) e(n\alpha_1) - e(-n\alpha_2) 
+\chi(-1)e(n\alpha_2)) + O(\sqrt{q}). \tag{7.1}
$$ 
There arise three cases: both $\alpha_1$ and $\alpha_2$ 
lie on minor arcs, exactly one of $\alpha_1$ and $\alpha_2$ 
lies on a major arc, and both $\alpha_1$ and $\alpha_2$ 
lie on major arcs.  In the first case we obtain 
from Lemma 6.5 that the above is 
$\ll \sqrt{q} (\log \log q)^{\frac 67+o(1)}$.  We 
examine the third case in detail, and omit  
the second case which is similar and simpler.  
Suppose (for $j=1$, $2$) that $|\alpha_j -b_j/r_j| \le 1/(r_j S)$ 
where $r_j \le s$, and $(b_j,r_j)=1$.  Set $N_j 
=\min (q, 1/|r_j \alpha_j -b_j|)$.  Using Lemma 6.5 and 
Lemma 4.4 we see that (7.1) equals, up to an error 
$O(\sqrt{q}(\log \log q)^{\frac 67+o(1)})$,
$$
-(1-\chi(-1))
\frac{\tau(\chi)}{2\pi i} \Big( \lambda_1 \sum\Sb v\le N_1 \\ 
v\in {\Cal S}((\log q)^2) \\ (v,r_1r_2)=1\endSb 
\frac{\cbar(v)}{v} - 
\lambda_2 \sum\Sb  v\le N_2 \\ 
v\in {\Cal S}((\log q)^2) \\ (v,r_1r_2)=1\endSb \frac{\cbar(v)}{v} 
\Big), \tag{7.2} 
$$
where 
$$
\lambda_j = \frac{1}{\phi(r_j)} \prod\Sb p^a \parallel r_j \\ a\ge 1 
\endSb (\cbar(p^a)-\cbar(p^{a-1}) ) \prod_{p|r_1r_2} \Big( 1-
\frac{\cbar(p)}{p}\Big)^{-1}.
$$
It is easy to see that (7.2) is bounded in magnitude by
$$
\align
&\frac{\sqrt{q}}{\pi} 
\max ( |\lambda_1|, |\lambda_2|, |\lambda_1-\lambda_2|) 
\sum\Sb v\in {\Cal S}((\log q)^2) \\ (v,r_1r_2)=1\endSb 
\frac 1v 
\\
=
&\frac{\sqrt{q}}{\pi} 
\max ( |\lambda_1|, |\lambda_2|, |\lambda_1-\lambda_2|) 
\prod_{p|r_1r_2} \Big(1-\frac 1p\Big) (2e^{\gamma}+o(1)) \log \log q.
\\
\endalign
$$
Thus Theorem 6 would follow if 
$$
\max(|\lambda_1|, |\lambda_2| ,
|\lambda_1-\lambda_2|) \phi(r_1r_2)/(r_1r_2) \le 1.  \tag{7.3}
$$

A simple optimization gives that 
$$
|\lambda_j| \frac{\phi(r_1r_2)}{r_1r_2} 
\le \frac{1}{r_j} \prod\Sb p^a \parallel r_j \\ a\ge 1\endSb 
\Big|\frac{\cbar(p^a)-\cbar(p^a-1)}{1-\cbar(p)/p}\Big| 
\le \frac{1}{r_j} \prod_{p|r_j}\frac{2p}{p+1} (\le 1).
$$
This estimate immediately gives (7.3) in all but the 
following two cases: one of $r_1$ or $r_2$ equals $1$,  or 
one of $r_1$ or $r_2$ equals $2$ and the other equals $3$.  
In the second case we see that 
$$
|\lambda_1-\lambda_2| \frac{\phi(6)}{6} 
= \Big|\frac{2\cbar(2)-\cbar(3)-1}{(2-\cbar(2))(3-\cbar(3))}\Big| 
\le \frac 23,
$$
since this is maximized at $\cbar(2)=-1$ and $\cbar(3)=1$ and 
so (7.3) holds.  Finally we have the case when one of 
$r_1$ or $r_2$ is $1$ and the other equals $r$ say.  Here 
we must show that 
$$
\Big| 1- \frac{1}{\phi(r)} \prod\Sb p^a \parallel 
r \\ a\ge 1\endSb (\cbar(p^a)-\cbar(p^{a-1})) \Big| 
\prod_{p|r} \Big|\frac{p-1}{p-\cbar(p)}\Big| 
\leq 1. \tag{7.4}
$$

If $r=p^a$ is a prime power then the LHS of (7.4) 
equals 
$$
\frac{1}{p^{a-1}}\Big|\Big(1-\frac 1p\Big) \frac{p^a-\cbar(p)^a}{p-\cbar(p)} 
+ \frac{\cbar(p)^{a-1}}{p}\Big| 
\le \frac{1}{p^{a-1}} \Big( \Big(1-\frac 1p\Big) (p^{a-1}+\ldots+1) 
+ \frac{1}{p} \Big) 
=1,
$$
and so (7.4) holds.  Now suppose 
that $r$ has at least two distinct prime factors.
For any non-negative $a_1,\dots,a_k$ with $k\geq 2$,
we have that
$$
(1+a_1\dots a_k)^2\leq (1+a_1^2)(1+(a_2\dots a_k)^2)\leq 
\prod_{i=1}^k (1+a_i^2).
$$
Therefore
$$
\align
\Big| 1 - \frac{1}{\phi(r)}\prod\Sb p^a\parallel r\\a\ge 1\endSb 
(\cbar(p^a)-\cbar(p^{a-1}))\Big|^2 &\le
\Big( 1+ \prod\Sb p^a\parallel r\\a\ge 1\endSb  
\frac{|\cbar(p^a)-\cbar(p^{a-1})|}{p^{a-1}(p-1)} \Big)^2
\le \prod_{p|r} \Big( 1+  \Big| \frac{1-\cbar(p)}{p-1}\Big|^2 \Big)\\
&\le \prod_{p|r} \Big| 1+  \frac{1-\cbar(p)}{p-1}\Big|^2
\endalign
$$
as desired, since  $z=(1-\cbar(p))/(p-1)$ is a complex number with
non-negative real part so that  $1+|z|^2\le |1+z|^2$.

\head 8. Paley's bound in all directions: Proof of 
Theorem 3 \endhead

\noindent  Bateman and Chowla
[1] showed that
$$
\frac{1}{q} \sum_{N\leq q} \Big| \sum_{n\leq N} \chi(n) -
\frac{\tau (\chi)}{i\pi} \frac{(1-\chi(-1))}{2}
L(1,\overline{\chi}) \Big|^2 =  \frac{q}{12} \prod_{p|q} \left(
1 - \frac{1}{p^2} \right) .
$$
If $\chi(-1)=-1$ we deduce that
$$
\sum_{n\leq N} \chi(n) = \frac{\tau (\chi)}{i\pi}
(L(1,\overline{\chi}) +O(\log\log\log q))  \tag{8.1} 
$$
for ``almost all'' $N\leq q$.  We now show that 
for most characters $\chi$, $L(1,\cbar)$ may 
be approximated by a short Euler product.  
Throughout this section we let $y:=\log q/\log \log q$.

\proclaim{Proposition 8.1} For any large prime $q$
$$
L(1,\chi)=\prod_{p\leq y}  \Big( 1 - \frac{\chi(p)}{p}
\Big)^{-1} \Big( 1+O\Big(\frac{\log \log \log q}{\log \log q}\Big) \Big),
$$
for all but at most $q^{1-1/(4\log\log q)}$ characters $\chi \pmod q$.
\endproclaim
\demo{Proof} An immediate consequence of
Proposition 2.2 of [8] is that
$$
L(1,\chi)=\prod_{p\leq (\log q)^3}  \left( 1 - \frac{\chi(p)}{p}
\right)^{-1} \Big(1 + O\Big(\frac{1}{\log \log q}\Big)\Big),
$$
for all but at most $q^{3/4}$ characters $\chi \pmod q$.  Consider 
$$
\frac{1}{\phi(q)} \sum_{\chi \pmod q}  \Big| \sum_{\log q <
p\le (\log q)^3} \frac{\chi(p)}{p} \Big|^{2k} = \sum_{m\equiv n
\pmod q} \frac{a_k(m)a_k(n)} {mn}
$$
where $a_k(n)$ is the number of ways of writing $n=p_1\dots p_k$
where each $p_i$ is a prime in $(\log q,(\log q)^3]$.  We 
choose $k=[\log q/(4\log \log q)]$ so that $a_k(n)=0$ for $n>q$ 
and so the congruence $m\equiv n\pmod q$ implies that $m=n$.  Since 
$a_k(n)\le k!$ it follows that 
$$
\sum_{n} \frac{a_k(n)^2} {n^2} \leq k! \sum_{n} \frac{a_k(n)}{n^2}
= k! \Big( \sum_{\log q < p \le (\log q)^3} \frac{1}{p^2}
\Big)^{k} 
\le \Big(\frac{1}{\log \log q}\Big)^{2k}.
$$
We deduce that there are fewer than $\phi(q)e^{-2k}$ characters $\chi$ 
with $|\sum_{\log q < p \le (\log q)^3} \chi(p)/p| \ge e/\log \log q$. 
Since $\sum_{y <p\le \log q} \frac{\chi(p)}{p} \ll \log \log \log q/\log 
\log q$ trivially,  the Proposition follows.
\enddemo

\proclaim{Proposition 8.2} Given a prime $q$ and an angle 
$\theta \in (-\pi, \pi]$,  there are at least 
$q^{1-C_0/(\log\log q)^2}$ characters $\chi \pmod q$
with $\chi(-1)=-1$ such that 
$$
\frac{\tau(\chi)}{i \sqrt{q}} \prod_{p\le y} 
\Big(1-\frac{\cbar(p)}{p}\Big)^{-1} 
= e^{i\theta} (e^{\gamma} \log \log q) 
+ O((\log \log q)^{1/2}).
$$
\endproclaim

\demo{Proof of Theorem 3}  Theorem 3 follows upon combining (8.1) 
with Propositions 8.1 and 8.2.  
\enddemo

To prove Proposition 8.2 we require the following 
consequence of P. Deligne's celebrated bound on 
hyper-Kloosterman sums. 

\proclaim{Lemma 8.3} We have 
$$ 
\frac{2}{\phi(q)} \Big|\sum\Sb \chi\pmod q\\ \chi(-1)=-1\endSb 
\chi(a) \tau(\chi)^n \Big|
\le 2n q^{(n-1)/2}.
$$
\endproclaim 
\demo{Proof}  Using the definition of $\tau(\chi)$ and 
the orthogonality relation for characters we 
see that 
$$
\frac{2}{\phi(q)}\sum\Sb \chi\pmod q\\ \chi(-1)=-1\endSb 
\chi(a) \tau(\chi)^n  = \text{Kl}_n(\overline{a},q) - 
\text{Kl}_n(-\overline{a},q),
$$
where  
$$
\text{Kl}_n(b,q) = \sum\Sb x_1, \ldots, x_n \pmod q \\ x_1\cdots x_n \equiv b
\pmod q \endSb e\Big(\frac{x_1+\ldots+x_n}{q}\Big).
$$
In (7.1.3) of [5] Deligne gives the bound (for $(b,q)=1$)
$$
|\text{Kl}_n(b,q)| \le n q^{(n-1)/2},
$$
from which the Lemma follows.

\enddemo 

\demo{Proof of Proposition 8.2} Set $R:=\prod_{p\le y}(1-1/p)^{-1} 
= e^{\gamma}\log \log q+ O(\log \log \log q)$ 
and consider for a natural number $k$
$$
\frac{2}{\phi(q)} \sum\Sb\chi \pmod q\\ \chi(-1)=-1 \endSb 
\Big|  \frac{\tau(\chi)}{i\sqrt{q}} \prod_{p\le 
y}\Big(1-\frac{\cbar(p)}{p}\Big)^{-1} 
+ Re^{i\theta}\Big|^{2k}. \tag{8.2}
$$
Expanding using the binomial theorem this 
equals 
$$
\align
&\sum_{0\le j, \ell \le k} \binom{k}{j}\binom{k}{\ell} 
R^{2k-j-\ell} e^{i\theta(\ell-j) } \sum\Sb m, n\in {\Cal S}(y) \endSb 
\frac{d_j(m)}{m}\frac{d_\ell(n)}{n} \\
&\hskip 1.5 in \times 
\frac{2}{\phi(q)} \sum\Sb \chi\pmod q\\ \chi(-1)=-1\endSb 
\cbar(m)\chi(n) \Big(\frac{\tau(\chi)}{i\sqrt{q}}\Big)^j 
\overline{\Big(\frac{\tau(\chi)}{i\sqrt{q}}\Big)}^{\ell}.
\tag{8.3}
\\
\endalign
$$

Using Lemma 8.3 we see that the terms $j\neq \ell$ above 
contribute an amount bounded in magnitude by
$$
\frac{2k}{\sqrt{q}} \sum_{0\le j,\ell \le k} 
\binom{k}{j}\binom{k}{\ell} R^{2k-j-\ell} R^j R^{\ell} 
= \frac{2k}{\sqrt{q}} 2^{2k} R^{2k}. \tag{8.4} 
$$
Now we focus on the terms $j=\ell$ in (8.3) which give, 
by the orthogonality relation for characters, 
$$ 
\sum_{0\le j\le k} \binom{k}{j}^2 R^{2k-2j} 
\sum\Sb m\equiv \pm n \pmod q \\ m, n\in {\Cal S}(y)\endSb 
(\pm 1) \frac{d_j(m)}{m} \frac{d_{j}(n)}{n}. \tag{8.5}
$$
If $m\equiv \pm n \pmod q$ but $m\neq n$ then either 
$m$ or $n$ exceeds $q/2$.  Thus such terms contribute 
to the sum in (8.5) an amount 
$$
\align
&\le 4\sum_{m \in {\Cal S}(y)} \frac{d_j(m)}{m} \sum\Sb n\ge q/2 \\ 
n\in {\Cal S}(y)\endSb \frac{d_{j}(n)}{n} 
\le 4 R^j \Big(\frac{2}{q}\Big)^{1/\log \log q} 
\sum_{n\in {\Cal S}(y)} \frac{d_j(n)}{n^{1-1/\log \log q}} 
\\
&\ll C^j R^{2j} q^{-1/\log \log q},
\\
\endalign
$$
for some absolute constant $C >1$.  From this and (8.4) we 
conclude that our quantity (8.2) is 
$$
\sum_{0\le j\le k}\binom{k}{j}^2 R^{2k-2j} 
\sum_{n\in {\Cal S}(y)} \frac{d_j(n)^2}{n^2} 
+ O\Big( (4C)^k R^{2k} q^{-1/\log \log q}\Big). \tag{8.6}
$$

Note that 
$$
\sum_{n \in {\Cal S}(y)} \frac{d_j(n)^2}{n^2} 
= \prod_{p\le y} \sum_{\ell =0}^{\infty} 
\frac{d_j(p^\ell)^2}{p^{2\ell} }
= \prod_{p\le y} \int_{-1/2}^{1/2} 
\Big| 1-\frac{e(\theta)}{p}\Big|^{-2j} 
d\theta.
$$
Observe that $\int_{-1/2}^{1/2} |1-\frac{e(\theta)}{p}|^{-2j} d\theta \ge 1$ 
always, and that if $p\le j$ then it is 
$$
\ge \int_{-p/(2j)}^{p/(2j)} \Big|1-\frac{e(\theta)}{p}\Big|^{-2j} d\theta
\ge \frac{cp}{j} \Big(1-\frac{1}{p} \Big)^{-2j},
$$
for a suitable positive constant $c$.  It follows that 
if $2\le j\le y$ then 
$$
\sum_{n \in {\Cal S}(y)} \frac{d_j(n)^2}{n^2} 
\ge \prod_{p\le j} \Big(1-\frac{1}{p}\Big)^{-2j} 
\exp\Big(- \frac{Cj}{\log j}\Big), \tag{8.7}
$$
for some positive constant $C$.  

We now take $k= [c'y]$ for a suitably small 
constant $c'>0$, and consider only the contribution of 
$j=[k/2]$ in (8.6).  Using (8.7) we deduce easily that the main term in 
(8.6) exceeds $2^{2k} R^k \prod_{p\le k} (1-1/p)^{-k} \exp(-Ck/\log k)$ 
and that the error term there is substantially smaller.  We conclude 
that 
$$
\frac{2}{\phi(q)} \sum\Sb\chi\pmod q\\ \chi(-1)=-1\endSb 
\Big| \frac{\tau(\chi)}{i\sqrt{q}} \prod_{p\le y} \Big( 1-\frac{\cbar(p)}{p}
\Big)^{-1} + Re^{i\theta}\Big|^{2k} 
\ge (2R)^{2k} \exp\Big(-\frac{Ck}{\log k}\Big),
$$
for some positive absolute constant $C$.  
From this estimate we deduce that (for some absolute 
constant $C_0$) there are at least 
$q^{1-C_0/(\log \log q)^2}$ characters $\chi \pmod q$ with 
$\chi(-1)=-1$ such that 
$$
\Big| \frac{\tau(\chi)}{i\sqrt{q}} \prod_{p\le y} \Big( 1-\frac{\cbar(p)}{p}
\Big)^{-1} + Re^{i\theta}\Big| 
\ge 2R \Big(1-\frac{C'}{\log \log q}\Big).
$$
Now, if $|z|\le 1$ and $|1+z| \ge 2-\epsilon$ then we 
may check easily that $z=1+O(\sqrt{\epsilon})$.  The 
Proposition follows.  

\enddemo


\head 9. The constant in the Polya-Vinogradov theorem: 
Proof of Theorem 2.7  \endhead

\noindent  Let $\xi \pmod m$ denote the primitive character 
with conductor below $(\log q)^{\frac 13}$ such that 
${\Bbb D}(\chi,\xi;q)$ is a minimum.  We distinguish two 
cases depending on whether $m>1$ or $m=1$.  

We start with the easier first case. By Lemma 3.4 we 
know that ${\Bbb D}(\chi,1;q)^2 \ge (1-1/\sqrt{2}+o(1)) 
\log \log q$, and so by (6.9) we have that 
$L(1,\cbar)= o(\log q)$.  Thus by 
(2.1) and Lemma 6.5 we deduce that 
$M(\chi) +o(\sqrt{q}\log q)$ is 
$$
\le \frac{\sqrt{qm}}{\pi} 
\max_{m|r, N\le q} \Big| \frac{1}{\phi(r)} 
\prod\Sb p^{a} \parallel r/m \\ a\ge 1\endSb 
(\cbar(p^a)-\overline{\xi}(p)\cbar(p^{a-1}))
\sum\Sb n\le N   \endSb \frac{(\cbar\xi)(n)}{n}\Big|. \tag{9.1}
$$

Observe that $\cbar \xi$ is 
a non-trivial character to the modulus $[q,m]$ (this denotes 
the l.c.m. of $q$ and $m$).  Set $c_0=1/4$ if $[q,m]$ is 
cubefree and $c_0=1/3$ otherwise, and note that 
$c_0=c$ unless $m$ is divisible by a prime cube, and in 
any case $c_0 \le \frac 43 c$.  Burgess's results on 
character sums (see [2]) show that $\sum_{n\le x} (\cbar \xi)(n) =o(x)$ 
if $x > q^{c_0+\epsilon}$, from which it follows by 
partial summation that 
$$
\sum_{n\le N} \frac{(\cbar \xi)(n)}{n} = 
\sum_{n\le \min(q^{c_0},N)} \frac{(\cbar \xi)(n)}{n} + o(\log q).
$$
Using this and Lemma 4.4 in (9.1) we obtain that 
$M(\chi) +o(\sqrt{q}\log q)$ is 
$$
\le \frac{\sqrt{qm}}{\pi} 
\max_{m|r, N\le q^{c_0}} \Big| \frac{1}{\phi(r)} 
\prod\Sb p^{a} \parallel r/m \\ a\ge 1\endSb 
\frac{\cbar(p^a)-\overline{\xi}(p)\cbar(p^{a-1})}{1-(\cbar\xi)(p)/p} 
\sum\Sb n\le N  \\ (n,r)=1 \endSb \frac{(\cbar\xi)(n)}{n}\Big|. \tag{9.2}
$$  

The product above is bounded in magnitude by $\prod_{p|r/m} 2p/(p+1)$, 
while the sum is bounded by $\sim (\phi(r)/r) c_0 \log q$.  Thus 
(9.2) is bounded in magnitude by 
$$
\frac{\sqrt{q}}{\pi} 
\frac{\sqrt{m}}{r} \prod_{p|r/m} \frac{2p}{p+1} (c_0 \log q).
$$
If $c_0\neq c$ then $m$ must be at least $8$ and the above 
bound beats the estimates claimed in the Theorem.  If 
$1<m<8$ then $c_0=c$ and the bound above 
suffices in all cases except for $m=r=3$ (and $\xi
=\fracwithdelims(){\cdot}{3}$).  In 
this final case we have that the quantity in (9.2) is bounded 
in magnitude by 
$$
\frac{\sqrt{q}}{\pi} \frac{\sqrt{3}}{2} \max_{N\le q^{c}}
\Big|\sum\Sb n\le N 
\endSb \frac{(\cbar \xi)(n)}{n}\Big|.
$$
Applying Theorem 1 of [9] we may see that
$$
|L(1,\cbar\xi)| = \Big| \sum_{n\le q^c} \frac{(\cbar\xi)(n)}{n}\Big| 
+o(\log q) \le \Big(\frac{34}{35} +o(1)\Big) \frac{2}{3} 
(c\log q),
$$
where the $\frac 23$ accounts for the fact that $(\cbar\xi)(3)=0$.  
It follows that for $N\le q^{c}$ 
$$
\align
\Big|\sum\Sb n\le N
\endSb \frac{(\cbar \xi)(n)}{n}\Big| &\le 
\left( \frac23+o(1)\right) \min\Big(   \log N,
\frac{34}{35}   (c\log q) +
\log \frac{q^c}{N}
\Big) \\
&\le (1+o(1)) \frac{69}{70}\cdot \frac{2}{3} (c\log q),\\
\endalign
$$
which completes the proof of the Theorem when $m>1$.  

Now consider the case $m=1$.  Here (2.1), Burgess's estimate, 
and Lemma 6.5 give that $M(\chi)+o(\sqrt{q}\log q)$ is 
$$
\le 
\frac{\sqrt{q}}{\pi} 
\max_{r, N\le q^c} \Big| \sum_{n\le q^c} \frac{\cbar(n)}{n} 
- \frac{1}{\phi(r)} \prod\Sb p^a \parallel r\endSb 
(\cbar(p^a)-\cbar(p^{a-1})) \sum\Sb n\le N \endSb 
\frac{\cbar(n)}{n}\Big|.
$$
The estimate claimed in the Theorem now follows 
from Lemma 4.4 and (7.4).

\enddocument